\documentclass{elsarticle}

\usepackage{lineno,hyperref}
\modulolinenumbers[5]

\journal{Applied Numerical Mathematics}

\usepackage{latexsym}
\usepackage{amsmath}
\usepackage{amssymb}
\usepackage{amsfonts}
\usepackage{verbatim}
\usepackage[latin1]{inputenc}
\usepackage{mathrsfs}
\usepackage{amsfonts}
\usepackage{graphicx}
\usepackage{colortbl,dcolumn}
\usepackage{amsmath}
\usepackage{psfrag}
\usepackage{booktabs}

\newtheorem{tm}{Theorem}[section]
\newtheorem{rk}{Remark}[section]

\newtheorem{prop}{Proposition}[section]
\newtheorem{lm}{Lemma}[section]

\newcommand{\bi} {\mathbf i}
\newcommand{\bk} {\mathbf k}

\newcommand{\PP}{\mathbb P}

\newcommand{\R}{\mathbb R}
\newcommand{\HH}{\mathbb H}

\newcommand{\Z}{\mathbb Z}
\newcommand{\CC}{\mathcal C}
\newcommand{\OO}{\mathcal O}

\begin{document}

\begin{frontmatter}

\title{Optimal Error Estimate of Conservative Local Discontinuous Galerkin Method for Nonlinear Schr\"odinger Equation}

\author[1]{Jialin Hong}
\ead{hjl@lsec.cc.ac.cn}

\author[2]{Lihai Ji}
\ead{jilihai@lsec.cc.ac.cn}

\author[3]{Zhihui Liu\corref{cor}}
\ead{liuzhihui@lsec.cc.ac.cn}

\cortext[cor]{Corresponding author.}

\address[1]{LSEC, ICMSEC, Academy of Mathematics and Systems Science, Chinese Academy of Sciences, Beijing 100190, China; School of Mathematical Sciences, University of Chinese Academy of Sciences, Beijing 100049, China}

\address[2]{Institute of Applied Physics and Computational Mathematics, Beijing 100094, China}

\address[3]{Academy of Mathematics and Systems Science, Chinese Academy of Sciences, Beijing 100049, China}

\begin{abstract}
In this paper, we propose a conservative local discontinuous Galerkin method for a one-dimensional nonlinear Schr\"odinger equation.
By using special generalized alternating numerical fluxes, we establish the optimal rate of convergence $\mathcal O(h^{\bk+1})$, with polynomial of degree $\bk$ and grid size $h$.
Meanwhile, we show that this method preserves the charge conservation law.
Numerical experiments verify our theoretical result.
\end{abstract}

\begin{keyword}
nonlinear Schr\"odinger equation \sep
optimal error estimates \sep
charge conservation law \sep
local discontinuous Galerkin method \sep
generalized alternating numerical flux
\MSC[2010]
Primary: 65M60 \sep
Secondary: 35Q55
\end{keyword}

\end{frontmatter}



\section{Introduction}
In this paper, we present a local discontinuous Galerkin (LDG) method with alternative numerical fluxes for focusing or defocusing nonlinear Schr\"odinger (NLS) equation
\begin{align}\label{nls}\tag{NLS}
\bi u_t+u_{xx}+\lambda |u|^2 u=0,
\quad (t,x)\in (0,T]\times \OO,
\end{align}
with initial datum $u(0,x)=u_0(x)$,
$x\in \OO:=(0,1)$, $\lambda\in \R$.
We will mainly focus on the periodic boundary condition $u(t,0)=u(t,1)$, $t\in [0,T]$.
It is well-known that Eq. \eqref{nls} possesses the charge conservation law, i.e.,
\begin{align*}
\|u(t)\|^2  =\|u_0\|^2, \quad \forall\ t\in (0,T].
\end{align*}
Our main observation is that the proposed LDG method preserves the charge and possesses an optimal convergence rate.
We call it a conservative local discontinuous Galerkin (CLDG) method.

The discontinuous Galerkin (DG) method is a class of finite element methods using discontinuous, piecewise
polynomials as the solution and the test spaces in the spatial direction. For a detailed description of the method as well as its implementation and applications, we refer the readers to the review paper \cite{CS01}. The LDG method is an
extension of the DG method aimed at solving partial differential equations (PDEs) containing higher than first order spatial derivatives. The idea of the LDG method is to rewrite the equations with higher order derivatives into a first order system, then apply the DG method
on the system. The design of the numerical fluxes is the key ingredient to ensure
stability. The LDG techniques have been developed for various high order PDEs, including convection diffusion equations \cite{CS98} and nonlinear one-dimensional and two-dimensional KdV-type equations \cite{XS10,YS02}. More details about the LDG methods for
high order time-dependent PDEs can be found in the review paper \cite{XS10}.

Since the basis functions can be completely discontinuous, the LDG methods have certain flexibilities and advantages.
It can be easily designed for any order of accuracy. In fact, the order of accuracy can
be locally determined in each cell, thus for efficient $h$-$p$ adaptivity. It is easy to handle complicated geometry and boundary conditions. It can be used
on arbitrary triangulations, even those with hanging nodes. It is extremely local in data communications. The evolution of the solution in each cell needs to communicate only with its immediate neighborhoods, regardless of the order of accuracy. The methods have excellent parallel efficiency. Finally, there is provable cell entropy inequality and $\mathbb L^2$ stability, for arbitrary scalar equations in any spatial dimension and any triangulation, for any order of accuracy, without limiters.

Some recent attempts have been made to apply the DG discretization to solve the Schr\"{o}dinger equation, see \cite{LCZ05, XS05, ZYF12, ZYZ12} and references therein. In \cite{XS05}, Xu and Shu developed an LDG method to solve the generalized NLS equation.
For linear Schr\"{o}dinger equation, they obtained an error estimate of order $k+\frac12$ for polynomials of degree $\bk$. In \cite{LCZ05}, Lu, Cai and Zhang presented an LDG method for solving one-dimensional linear Schr\"{o}dinger equation so that the mass is preserved numerically. Zhang, Yu and Feng presented a mass preserving direct discontinuous Galerkin method in \cite{ZYF12} for the one-dimensional coupled NLS equations, and in \cite{ZYZ12} for both one and two dimensional NLS equations. Particularly, in \cite{ZYZ12} the conservation property is verified, and further validated by some long time simulation results.

Compared with the status of optimal $\mathbb L^2$-error estimates for LDG methods solving
time-dependent diffusive PDEs, for example, the convection diffusion equations \cite{CS98,Shu09}, optimal $\mathbb L^2$-error estimates for LDG methods solving high order time-dependent wave equations
are much more elusive. The main technical difficulty is the lack of coercivity and hence
the control on the auxiliary variables in the LDG method which are approximations
to the derivatives of the solution and the lack of control on the interface boundary
terms. When these issues are not addressed carefully, optimal $\mathbb L^2$-error estimates
could not be obtained. In \cite{XS05,YS02}, a priori $\mathbb L^2$-error estimates with suboptimal order
$k+\frac12$ for the LDG method with $\mathbb P^{k}$ elements for the linearized KdV equations and
the linearized Schr\"{o}dinger equation in one spatial dimension were obtained. For high order linear wave equations, \cite{XS12} proposed a general approach for proving optimal error estimates $k+1$ by utilizing the LDG method and its time derivatives with different test functions and fully making use of the so-called Gauss-Radau projections. In \cite{XCS13}, the authors developed an energy conserving LDG method for solving the second order linear wave equation and showed an optimal error estimate.
In \cite{CMZ16}, the authors consider the LDG method for solving the linear convection-diffusion equations and obtain directly the optimal $\mathbb L^2$-norm error estimate in a uniform
framework.
Recently, \cite{MSW16} presented an optimal $\mathbb L^2$-error estimate of the LDG method based on upwind-biased numerical fluxes for linear hyperbolic problems.

The aim of this paper is to obtain the optimal rate of convergence order $k+1$ for the CLDG method with a generalized numerical flux and a special projections
on the auxiliary variables.
The optimal error estimates hold not only for the solution itself but also for the auxiliary variables in the CLDG method approximating the various order derivatives of the solution. To our best knowledge, this is the first successful optimal $\mathbb L^2$-error estimates of the CLDG methods for such high order equations when not purely upwind numerical fluxes are considered.
We also note that the arguments in the present paper can be adapted to general NLS equation
\begin{align}\label{nls0}
\bi u_t+u_{xx}+f(|u|^2) u=0 \quad\text{in}\quad (0,T]\times \OO
\end{align}
with sufficiently smooth, real-valued function $f$.

The paper is organized as follows.
In Section \ref{sec-2}, we present the CLDG method for Eq. \eqref{nls} and their well-posedness and a priori estimations.
In Section \ref{sec-3}, we show the CLDG method possesses the charge conservation law and optimal convergence rate results as well as generalization of Eq. \eqref{nls0}.
Numerical experiments confirming the optimality of our theoretical results are given in Section \ref{sec-4}.
Concluding remarks are given in Section \ref{sec-5}.

\section{CLDG Method for NLS Equation}
\label{sec-2}

In this section we introduce notations and definitions to be used later in the paper and propose a CLDG method for Eq. \eqref{nls}.

\subsection{Basic Notations}

Define $\Z_N:=\{1,2,\cdots,N\}$ for an integer $N\geq 1$.
We denote by $\mathcal I_h$ a resellation of $\OO=(0,1)$ devided into $N$ cells $\OO_j=(x_{j-\frac12},x_{j+\frac12})$, and denote by $x_j=\frac12(x_{j-\frac12}+x_{j+\frac12})$ its center, $j\in \Z_N$.
Let $h_j=x_{j+\frac12}-x_{j-\frac12}$
and $h=\max_{j\in \Z_N} h_j$.
Assume that the mesh is quasi-uniform in the sense that there exists a positive constant $\gamma$ such that $\gamma h\le h_j$ for any $j\in \Z_N$.

For an positive integer $\bk$, we define a finite-element space consisting of piecewise polynomials
\begin{align*}
V_h^\bk:=\{v\in \mathbb L^2(\OO):v\big|_{\OO_j}\in \PP^\bk(\OO_j),\ j\in \Z_N\},
\end{align*}
where $\mathbb{P}^\bk(\OO_j)$ denotes the space of polynomials of the degree up to $\bk$ in each cell $\OO_j$. Note that functions in $V_h^\bk$ are allowed to be discontinuous across element interfaces.
The solution of the numerical method is denoted by $u_h$ which belongs to $V_h^\bk$. We denote by $(u_h)^-_{j+\frac12}$ and $(u_h)^+_{j+\frac12}$ the left and right limits of $u_h$ at $x_{j+\frac12}$, respectively.

\subsection{CLDG Method}

In order to construct the CLDG method, we rewrite \eqref{nls} as the first-order system
\begin{align}\label{ldg}
\begin{split}
\bi u_t+v_x+\lambda |u|^2 u
&=0,\\
v-u_x
&=0.
\end{split}
\end{align}
The LDG method for solving \eqref{ldg} is defined as follows: find $u_h,v_h\in V_h^\bk$ such that for all test functions $\alpha_h,\beta_h\in V_h^\bk$ and all $j\in \Z_N$,
\begin{align}\label{ldg0}
\begin{split}
\bi  \int_{\OO_j}(u_h)_t \alpha_h {\rm d}x-\int_{\OO_j}v_h (\alpha_h)_x {\rm d}x
+(\widehat{v}_h \alpha_h^-)_{j+\frac12}-(\widehat{v}_h \alpha_h^+)_{j-\frac12}\\
+\lambda \int_{\OO_j} |u_h|^2 u_h \alpha_h {\rm d}x=0&,\\
\int_{\OO_j}v_h \beta_h {\rm d}x+\int_{\OO_j} u_h (\beta_h)_x {\rm d}x-(\widehat{u}_h \beta_h^-)_{j+\frac12}+(\widehat{u}_h \beta_h^+)_{j-\frac12}=0&.
\end{split}
\end{align}

In this paper, instead of using the purely upwind flux, we adopt a generalized alternating numerical flux.
To be more specific, we choose
\begin{align}\label{flux}
\begin{split}
\widehat{u}_h=\theta u_h^-+(1-\theta)u_h^+\quad \text{at}\quad x_{j+\frac12},\quad j=\{0\}\cup \Z_N,\\
\widehat{v}_h=(1-\theta) v_h^-+\theta v_h^+\quad \text{at}\quad x_{j+\frac12},\quad j=\{0\}\cup \Z_N,
\end{split}
\end{align}
where $\theta\in [0,1]$.

We decompose the complex function $u(t,x)$ into real and imaginary parts:
\begin{align*}
u(t,x)=r(t,x)+\bi s(t,x),\quad (t,x)\in [0,T]\times \OO,
\end{align*}
where $r$ and $s$ being real-valued functions. Under the new notation, Eq. \eqref{nls} can be written as
\begin{align*}
r_t+s_{xx}+\lambda (r^2+s^2) s=0,\\
s_t-r_{xx}+\lambda (r^2+s^2) r=0,
\end{align*}
which is equivalent to the first-order system
\begin{align*}
\begin{split}
p-s_x&=0,\\
r_t+p_x+\lambda (r^2+s^2) s&=0,\\
q-r_x&=0,\\
s_t-q_x+\lambda (r^2+s^2) r&=0.
\end{split}
\end{align*}
The LDG method \eqref{ldg0} is equivalent to find $r_h,p_h,s_h,q_h\in V_h^\bk$ such that for any 
$\gamma_h,\omega_h,\alpha_h,\beta_h\in V_h^\bk$,
\begin{align}\label{ldg1}\tag{CLDG1}
\begin{split}
\int_{\OO_j}(r_h)_t \gamma_h {\rm d}x
-\int_{\OO_j}p_h (\gamma_h)_x {\rm d}x
+\Big[(\widehat{p}_h \gamma_h^-)_{j+\frac12}
-(\widehat{p}_h \gamma_h^+)_{j-\frac12}\Big] \\
+\lambda \int_{\OO_j} (r_h^2+s_h^2) s_h \gamma_h {\rm d}x=0&,\\
\int_{\OO_j}p_h \omega_h {\rm d}x+\int_{\OO_j}s_h (\omega_h)_x {\rm d}x
-\Big[(\widehat{s}_h \omega_h^-)_{j+\frac12}
-(\widehat{s}_h \omega_h^+)_{j-\frac12} \Big] =0&,\\
\int_{\OO_j}(s_h)_t \alpha_h {\rm d}x
+\int_{\OO_j}q_h (\alpha_h)_x {\rm d}x
-\Big[(\widehat{q}_h \alpha_h^-)_{j+\frac12}
-(\widehat{q}_h \alpha_h^+)_{j-\frac12}\Big] \\
-\lambda \int_{\OO_j} (r_h^2+s_h^2) r_h \alpha_h {\rm d}x=0&,\\
\int_{\OO_j}q_h \beta_h {\rm d}x+\int_{\OO_j}r_h (\beta_h)_x {\rm d}x
-\Big[(\widehat{r}_h \beta_h^-)_{j+\frac12}
+(\widehat{r}_h \beta_h^+)_{j-\frac12}\Big]=0&,
\end{split}
\end{align}
and the numerical fluxes become
\begin{align}\label{flux1}\tag{CLDG2}
\begin{split}
\widehat{r}_h=\theta r_h^-+(1-\theta) r_h^+\quad \text{at}\quad x_{j+\frac12},\quad j=\{0\}\cup \Z_N,\\
\widehat{p}_h=(1-\theta) p_h^-+\theta p_h^+\quad \text{at}\quad x_{j+\frac12},\quad j=\{0\}\cup \Z_N,\\
\widehat{s}_h=\theta s_h^-+(1-\theta) s_h^+\quad \text{at}\quad x_{j+\frac12},\quad j=\{0\}\cup \Z_N,\\
\widehat{q}_h=(1-\theta) q_h^-+ \theta q_h^+\quad \text{at}\quad x_{j+\frac12},\quad j=\{0\}\cup \Z_N.
\end{split}
\end{align}

\subsection{Well-posedness and A Priori Estimations}

We denote by $\HH^\bk(\OO)$ the standard Sobolev space endowed with the norm
\begin{align*}
\|u\|_{\HH^\bk (\OO)}
:=\sum_{l=0}^\bk \|u^{(l)}(x)\|_{\mathbb L^2(\OO)}.
\end{align*}
To derive the optimal convergence rate of the CLDG method \eqref{ldg1}--\eqref{flux1}, we need the following a priori estimations for Eq. \eqref{nls}.

\begin{lm}\label{well}
Assume that $u_0\in \HH^{\bk+3}$. Eq. \eqref{nls} exists a unique solution $u$ such that
\begin{align*}
u\in \CC([0,T]; \HH^{\bk+3}(\OO))
\quad \text{and} \quad
u_t\in \CC([0,T]; \HH^{\bk+1}(\OO)).
\end{align*}
As a consequence, $u\in \CC([0,T]; \mathbb L^\infty(\OO))$.
\end{lm}

{\it Proof.}
It is known that $u\in \CC([0,T]; \HH^{\bk+3}(\OO))$ and $u_t\in \CC([0,T]; \HH^{\bk+1}(\OO))$ (see \cite{GV84}, Proposition 1.2).
By Sobolev embedding $\HH^1\hookrightarrow \mathbb L^\infty$, it follows that $u\in \CC([0,T]; \mathbb L^\infty(\OO))$.
\quad \\

\section{Main Results}
\label{sec-3}

\subsection{Charge Conservation Law}

In this subsection, we present the charge conservation law of LDG method \eqref{ldg0}.
\begin{prop}
{\em There exist numerical entropy fluxes $\hat{\phi}_{j+\frac12}$ such that the solution to the method \eqref{ldg0}--\eqref{flux} satisfies
\begin{align}\label{ent}
&\frac d{dt}\left[ \int_{\OO_j}|u_h|^2 {\rm d}x\right]
+\hat{\phi}_{j+\frac12}-\hat{\phi}_{j-\frac12}=0.
\end{align}}
\end{prop}

{\it Proof.}
First, we take the complex conjugate for every term in \eqref{ldg0} and obtain
\begin{align}\label{ldg2}
\begin{split}
-\bi  \int_{\OO_j}(u_h^*)_t \alpha_h^* {\rm d}x
-\int_{\OO_j}v_h^* (\alpha_h^*)_x {\rm d}x
+(\widehat{v}_h^* \alpha_h^{*-})_{j+\frac12}
-(\widehat{v}_h^* \alpha_h^{*+})_{j-\frac12} \\
+\lambda \int_{\OO_j} |u_h|^2 u_h^* \alpha_h^*{\rm d}x=0&,\\
\int_{\OO_j}v_h^* \beta_h^* {\rm d}x+\int_{\OO_j} u_h^* (\beta_h^*)_x {\rm d}x
-(\widehat{u}_h^* \beta_h^{*-})_{j+\frac12}
+(\widehat{u}_h^* \beta_h^{*+})_{j-\frac12}=0&,
\end{split}
\end{align}
where $u_h^*$ denotes the complex conjugate of $u_h$. Since \eqref{ldg0} and \eqref{ldg2} hold for any test functions in $V_h^\bk$, we choose
$\alpha_h=u_h^*,~~\beta_h=v_h^*$.
With these choices of test functions, it follows that
\begin{align}\label{ldg22}
\begin{split}
\bi  \int_{\OO_j}(u_h)_t u_h^* {\rm d}x
-\int_{\OO_j}v_h (u_h^*)_x {\rm d}x
+(\widehat{v}_h u_h^{*-})_{j+\frac12}
-(\widehat{v}_h u_h^{*+})_{j-\frac12} \\
+\lambda \int_{\OO_j} |u_h|^4 {\rm d}x=0&,\\
\int_{\OO_j}v_h v_h^* {\rm d}x+\int_{\OO_j} u_h (v_h^*)_x {\rm d}x
-(\widehat{u}_h v_h^{*-})_{j+\frac12}
+(\widehat{u}_h v_h^{*+})_{j-\frac12}=0&.
\end{split}
\end{align}
and
\begin{align}\label{ldg21}
\begin{split}
-\bi  \int_{\OO_j}(u_h^*)_t u_h {\rm d}x
-\int_{\OO_j}v_h^* (u_h)_x {\rm d}x
+(\widehat{v}_h^* u_h^-)_{j+\frac12}
-(\widehat{v}_h^* u_h^+)_{j-\frac12} \\
+\lambda \int_{\OO_j} |u_h|^4 {\rm d}x=0&,\\
\int_{\OO_j}v_h^* v_h {\rm d}x+\int_{\OO_j} u_h^* (v_h)_x {\rm d}x
-(\widehat{u}_h^* v_h^-)_{j+\frac12}
+(\widehat{u}_h^* v_h^+)_{j-\frac12}=0&.
\end{split}
\end{align}
Adding the two equalities in \eqref{ldg22}, we get
\begin{align*}
\begin{split}
0&=\bi  \int_{\OO_j}(u_h)_t u_h^* {\rm d}x
-\int_{\OO_j}v_h (u_h^*)_x {\rm d}x+\int_{\OO_j}v_h v_h^* {\rm d}x+\int_{\OO_j} u_h (v_h^*)_x {\rm d}x \\
&\quad +\lambda \int_{\OO_j} |u_h|^4 {\rm d}x
+(\widehat{v}_h u_h^{*-})_{j+\frac12}
-(\widehat{v}_h u_h^{*+})_{j-\frac12}
-(\widehat{u}_h v_h^{*-})_{j+\frac12}
+(\widehat{u}_h v_h^{*+})_{j-\frac12}.
\end{split}
\end{align*}
Similarly, for \eqref{ldg21}, it yields
\begin{align*}
\begin{split}
0 &=-\bi  \int_{\OO_j}(u_h^*)_t u_h {\rm d}x
-\int_{\OO_j}v_h^* (u_h)_x {\rm d}x+\int_{\OO_j}v_h^* v_h {\rm d}x+\int_{\OO_j} u_h^* (v_h)_x {\rm d}x \\
&\quad +\lambda \int_{\OO_j} |u_h|^4 {\rm d}x
+(\widehat{v}_h^* u_h^-)_{j+\frac12}
-(\widehat{v}_h^* u_h^+)_{j-\frac12}
-(\widehat{u}_h^* v_h^-)_{j+\frac12}
+(\widehat{u}_h^* v_h^+)_{j-\frac12}.
\end{split}
\end{align*}
Then taking the difference between the above two equalities leads to
\begin{align}\label{ldg21-22}
\begin{split}
0=&\bi  \int_{\OO_j}\left[(u_h)_t u_h^*+(u_h^*)_t u_h\right] {\rm d}x \\
&-\int_{\OO_j}\left[v_h (u_h^*)_x-v_h^* (u_h)_x-u_h (v_h^*)_x+u_h^* (v_h)_x\right] {\rm d}x \\
&-(\widehat{u}_h v_h^{*-})_{j+\frac12}
+(\widehat{u}_h^* v_h^-)_{j+\frac12}
-(\widehat{v}_h^* u_h^-)_{j+\frac12}
+(\widehat{v}_h u_h^{*-})_{j+\frac12} \\
&+(\widehat{u}_h v_h^{*+})_{j-\frac12}
-(\widehat{u}_h^* v_h^+)_{j-\frac12}
+(\widehat{v}_h^* u_h^+)_{j-\frac12}
-(\widehat{v}_h u_h^{*+})_{j-\frac12}.
\end{split}
\end{align}

For the second term, we have
\begin{align*}
&\int_{\OO_j}\left[v_h (u_h^*)_x-v_h^* (u_h)_x-u_h (v_h^*)_x+u_h^* (v_h)_x\right] {\rm d}x\\
&=(u_h^- v_h^{*-})_{j+\frac12}-(u_h^+ v_h^{*+})_{j-\frac12}-(v_h^- u_h^{*-})_{j+\frac12}+(v_h^+ u_h^{*+})_{j-\frac12}.
\end{align*}
Substituting it into Eq. \eqref{ldg21-22}, we obtain
\begin{align*}
&\frac d{dt} \left[ \int_{\OO_j}|u_h|^2 {\rm d}x \right]
+\hat{\phi}_{j+\frac12}-\hat{\phi}_{j-\frac12}=0,
\end{align*}
where the numerical entropy flux is given by
$$\hat{\phi}=2\text{Im}\left( \theta v_h^+u_h^{*-}+(1-\theta) v_h^- u_h^{*+}\right).$$
This completes the proof.
\quad \\

\begin{tm}
The solution to the LDG method \eqref{ldg0}--\eqref{flux} possesses the charge conservation law, i.e.,
\begin{align}\label{ccl}
\|u_h(t)\|^2  =\|u_0\|^2, \quad \forall\ t\in (0,T].
\end{align}
\end{tm}

{\it Proof.}
Summing up Eq. \eqref{ent} with $j$ over $\Z_N$ and using the periodic boundary condition, we have
\begin{align*}
&\frac d{dt} \left[ \int_{\OO} |u_h|^2 {\rm d}x \right]=0,
\end{align*}
from which we obtain \eqref{ccl}.
\quad \\

\begin{rk}
We call the LDG method \eqref{ldg0}--\eqref{flux} the CLDG method.
The charge conservation law trivially implies an $\mathbb L^2$-stability of the numerical solution.
\end{rk}

\subsection{Optimal Error Estimates}

In this subsection, we obtain the optimal error estimates for the approximations $r_h, s_h\in V_h^\bk$, which are given by the CLDG method \eqref{ldg1}--\eqref{flux1}.

\subsubsection{Projection and Interpolation Properties}\label{pro_p}

In what follows, we consider two special projections of a function $u$ with $k+1$ continuous derivatives into the space $V_h^\bk$. The special projections $\mathcal{P},~\mathcal{Q}$ are defined as follows. Given a function $u\in H^1(\mathcal{I}_h)$ and any subinterval $\OO_j$, it holds that
\begin{align}\label{pro1}
\begin{split}
&\int_{\OO_j} \left[\mathcal Pu(x)-u(x)\right] \omega {\rm d}x=0
\quad \forall\ \omega\in \mathbb{P}^{k-1}(\OO_j),\\
&\widehat{\mathcal Pu}_{j+\frac12}=\widehat{u}_{j+\frac12},
\quad j\in \Z_N
\end{split}
\end{align}
and
\begin{align}\label{pro2}
\begin{split}
&\int_{\OO_j} \left[\mathcal Qu(x)-u(x)\right] \omega {\rm d}x=0
\quad \forall\ \omega\in \mathbb{P}^\bk(\OO_j),\\
&\widehat{\mathcal Qu}_{j+\frac12}=\widehat{u}_{j+\frac12},
\quad j\in \Z_N.
\end{split}
\end{align}
Here and below, we denote $\hat{\omega}:=\theta\omega^-+(1-\theta)\omega^+$ for any $\omega\in H^{1}(\mathcal I_h)$.
In particular, when $\theta=1$, $\mathcal P$, $\mathcal Q$ are Gauss-Radau projection $\mathcal P^-$ and $\mathcal Q^-$, respectively.
It is well-known that (see e.g. \cite{Cia78}, Theorem 3.1.6) there holds for any
$j\in \Z_N$ that
\begin{align}\label{g-r}
\begin{split}
\|u-\mathcal P^- u\|_{\mathbb L^\infty(\OO_j)}
&\le Ch^{\bk+1} \|u\|_{\HH^{\bk+1} (\OO_j)},\\
\|u-\mathcal Q^- u\|_{\mathbb L^\infty(\OO_j)}
&\le Ch^{\bk+1} \|u\|_{\HH^{\bk+1} (\OO_j)}.
\end{split}
\end{align}

The projections mentioned above are shown in \cite[Lemma 2.6]{MSW16} to be well-defined.
Indeed, denote by $\mathcal P^-$ the Gauss-Radau projection and $E:=\mathcal P-\mathcal P^-$.
Since $\mathcal P^-$ is unique, the existence and uniqueness of $\mathcal P$ are equivalent to those of $E$.
\cite{MSW16} has proved that $E_j$, the restriction of $E$ to each $\OO_j$, can be represented as
\begin{align*}
E_j(x)=\sum_{l=0}^\bk \alpha_{j,l} P_{j,l}(x)
=\sum_{l=0}^\bk \alpha_{j,l} P_l(\xi),
\end{align*}
where $P_l(\xi)$ are the $l$-order Legendre polynomials and are orthogonal on $[-1,1]$ with $\xi=2(x-x_j)/h_j$ and on each element,
$P_{j,l}(x):=P_l(\xi)$ for $x\in \OO_j$ and the coefficients satisfy
\begin{align*}
\alpha_{j,l}=0,\quad l=0,1,\cdots,k-1;\ j=1,\cdots,N.
\end{align*}
Moreover, if we define $\eta_{j+1}:=(u-\mathcal P^- u)_{j+\frac12}^+$ for $j=\{0\}\cup \Z_{N-1}$ with $\eta_{N+1}=\eta_1$,
$\eta=(\eta_2,\eta_3,\cdots, \eta_{N+1})^T$ and
$\alpha_k=(\alpha_{1,k},\alpha_{2,k}, \cdots, \alpha_{N,k})^T$, then
\begin{align*}
A\alpha_k=(1-\theta)\eta,
\end{align*}
where $A=\text{circ}(\theta,(1-\theta)(-1)^\bk,0,\cdots,0)$ is an $N\times N$ circulant matrix.
The determinant of $A$ is
\begin{align} \label{det-A}
|A|=\theta^N (1-q^N)\quad \text{with}\quad
q=\frac{(1-\theta)(-1)^{\bk+1}}{\theta},
\end{align}
from which we conclude $A$ is always invertible for all $\bk$ and $N$ whenever $\theta\neq 1/2$.
This establishes existence and uniqueness of $\mathcal P$.

Moreover, \cite{CMZ16} obtains the following estimates of these projections.
The proof is similar to \cite{CMZ16}, Lemma 3.2, and we omit the details.

\begin{lm}\label{pro-est}
Assume that $u\in \HH^{\bk+1} (\mathcal{I}_h)$.
For any $\theta\neq 1/2$, there exists $C=C(\theta)$ which is independent of $h$ such that
\begin{align}\label{pro-est0}
\begin{split}
\|u-\mathcal P u\|_{\mathbb L^2(\OO)}
\le C h^{\bk+1} \|u\|_{\HH^{\bk+1} (\OO)},\\
\|u-\mathcal Q u\|_{\mathbb L^2(\OO)}
\le C h^{\bk+1} \|u\|_{\HH^{\bk+1} (\OO)}.
\end{split}
\end{align}
\end{lm}

\begin{rk}\label{rk-12}
When $\theta=1/2$, then by \eqref{det-A} we have
 $|A|=2^{-N} ( 1-(-1)^{N(\bk+1)})$, which shows that
$A$ is invertible if and only if $N$ is odd and $\bk$ is even.
In this case, $q=-1$ and $\beta_N=1$.
Then
\begin{align*}
|\alpha_{j,k}|
\le CN h^{\bk+1} \|u\|_{\HH^{\bk+1} (\OO)}
\le C h^\bk\|u\|_{\HH^{\bk+1} (\OO)},
\end{align*}
and thus
\begin{align*}
\|u-\mathcal P u\|_{\mathbb L^2(\OO)}
\le C h^\bk \|u\|_{\HH^{\bk+1} (\OO)},\\
\|u-\mathcal Q u\|_{\mathbb L^2(\OO)}
\le C h^\bk \|u\|_{\HH^{\bk+1} (\OO)}.
\end{align*}
Compared with the estimate \eqref{pro-est} for $\theta\neq 1/2$, the order of projection for $\theta=1/2$ reduces one.
However, it should be noted that in the numerical experiments we observe that the case $\theta=1/2$ also achieves the optimal error estimate.
\end{rk}

\subsubsection{Notations for the CLDG Discretization}

To facilitate the proof of the error estimate, we define the CLDG descretization operator $\mathcal{B}$, i.e., for each subinterval $\OO_j$,
\begin{align*}
\begin{split}
&\mathcal B_{j}(r,p,s,q;~\gamma,\omega,\alpha,\beta) \\
=&-\int_{\OO_j}p \gamma_x {\rm d}x+\int_{\OO_j} s\omega_x {\rm d}x+\int_{\OO_j}q \alpha_x {\rm d}x
+\int_{\OO_j} r \beta_x {\rm d}x \\
&+(\widehat{p} \gamma^-)_{j+\frac12}
-(\widehat{p} \gamma^+)_{j-\frac12}
-(\widehat{s} \omega^-)_{j+\frac12}
+(\widehat{s} \omega^+)_{j-\frac12}  \\
&-(\widehat{q} \alpha^-)_{j+\frac12}
+(\widehat{q} \alpha^+)_{j-\frac12}
-(\widehat{r} \beta^-)_{j+\frac12}
+(\widehat{r} \beta^+)_{j-\frac12},
\end{split}
\end{align*}
and
\begin{align*}
\mathcal B(r,p,s,q;~\gamma,\omega,\alpha,\beta)
:=\sum_{j} \mathcal B_j(r,p,s,q;~\gamma,\omega,\alpha,\beta).
\end{align*}

According to the periodic boundary condition and the definitions of the operator $\mathcal B$ and the projections $\mathcal P,\ \mathcal Q$, we have

\begin{lm}\label{lemma1}
{\em For any $\gamma_h,\omega_h,\alpha_h,\beta_h\in V_h^\bk$, it holds
\begin{align}\label{pro}
\mathcal B(r-\mathcal P r,p-\mathcal Q   p,s-\mathcal P s,q-\mathcal Q   q; ~\gamma_h,\omega_h,\alpha_h,\beta_h)=0,
\end{align}
where $\mathcal{P},~\mathcal{Q}$ are the projections defined in section \ref{pro_p}.}
\qed
\end{lm}

In addition, we define the LDG discretization operator $\mathcal{H}$ for the nonlinear term,  i.e., for each subinterval $\OO_j$,
\begin{align*}
\mathcal H_{j}(r,s;~\gamma,\alpha)
:=\lambda \int_{\OO_j} (r^2+s^2)r\alpha{\rm d}x
-\lambda \int_{\OO_j} (r^2+s^2)s \gamma{\rm d}x
\end{align*}
and
\begin{align*}
\mathcal H(r,s;~\gamma,\alpha)
:=\sum_{j} \mathcal H_j(r,s;~\gamma,\alpha).
\end{align*}

\subsubsection{Optimal Error Estimates}

With these notations and equalities \eqref{ldg1}--\eqref{flux1}, we have
\begin{align*}
\begin{split}
&\int_{\OO_j}(r_h)_t \gamma_h {\rm d}x+\int_{\OO_j}(s_h)_t \alpha_h {\rm d}x+\int_{\OO_j}p_h\omega_h {\rm d}x
+\int_{\OO_j}q_h\beta_h {\rm d}x\\
&\quad
+\mathcal B_j(r_h,p_h,s_h,q_h;~\gamma_h,\omega_h,\alpha_h,\beta_h)-\mathcal H_j(r_h,s_h;~\gamma_h,\alpha_h)=0.
\end{split}
\end{align*}

Notice that the CLDG method \eqref{ldg1}--\eqref{flux1} is also satisfied when the numerical solutions $r_h,p_h,s_h,q_h$ are replaced by the exact solutions $r$, $p=r_x$, $s$, $q=s_x$.
This yields the following cell error equation:
\begin{align*}
\begin{split}
&\int_{\OO_j}(r-r_h)_t \gamma_h {\rm d}x+\int_{\OO_j}(s-s_h)_t \alpha_h {\rm d}x+\int_{\OO_j}(p-p_h)\omega_h {\rm d}x\\
&+\int_{\OO_j}(q-q_h)\beta_h {\rm d}x+\mathcal{B}_j(r-r_h,p-p_h,s-s_h,q-q_h;~\gamma_h,\omega_h,\alpha_h,\beta_h)\\
&-\mathcal{H}_j(r,s;~\gamma_h,\alpha_h)+\mathcal{H}_j(r_h,s_h;~\gamma_h,\alpha_h)=0,
\quad \gamma_h,\omega_h,\alpha_h,\beta_h\in V_h^\bk.
\end{split}
\end{align*}
Summing over $j$, we get the error equation
\begin{align*}
  \begin{split}
  &\int_{\OO}(r-r_h)_t \gamma_h {\rm d}x+\int_{\OO}(s-s_h)_t \alpha_h {\rm d}x+\int_{\OO}(p-p_h)\omega_h {\rm d}x\\
  &+\int_{\OO}(q-q_h)\beta_h {\rm d}x+\mathcal{B}(r-r_h,p-p_h,s-s_h,q-q_h;~\gamma_h,\omega_h,\alpha_h,\beta_h)\\
  &=\mathcal{H}(r,s;~\gamma_h,\alpha_h)-\mathcal{H}(r_h,s_h;~\gamma_h,\alpha_h),\quad \gamma_h,\omega_h,\alpha_h,\beta_h\in V_h^\bk.
  \end{split}
\end{align*}

\begin{tm}\label{u-uh}
Assume that $u_0\in \HH^{\bk+3}$.
Let $u$ and $u_h$ be the solutions of Eq. \eqref{nls} and Eq. \eqref{ldg1}--\eqref{flux1}, respectively, with $\theta\neq 1/2$.
There exists a constant $C=C(T,u_0,\theta)$ such that
\begin{align}\label{u-uh0}
  \|u-u_h\|_{\mathbb L^2(\OO)}\le C(1+\|u_h\|_{\mathbb L^\infty(\OO)})h^{\bk+1}.
\end{align}
\end{tm}

{\it Proof.}
Denote
\begin{align*}
    &e_r=r-r_h=r-\mathcal{P}r+\mathcal{P}e_r,\\
    &e_s=s-s_h=s-\mathcal{P}s+\mathcal{P}e_s,\\
    &e_p=p-p_h=p-\mathcal{Q}p+\mathcal{Q}e_p,\\
    &e_q=q-q_h=q-\mathcal{Q}q+\mathcal{Q}e_q.
\end{align*}
Taking the test functions
\begin{equation*}
\gamma_h=\mathcal{P}r-r_h,~\omega_h=\mathcal{Q}p-p_h,~\alpha_h=\mathcal{P}s-s_h,~\beta_h=\mathcal{Q}q-q_h,
\end{equation*}
we obtain
\begin{align*}
\begin{split}
&\frac12\frac{d}{dt} \left[ \int_{\OO}(\mathcal P e_r)^2{\rm d}x\right]
+\frac12\frac{d}{dt} \left[\int_{\OO}(\mathcal P e_s)^2{\rm d}x\right]
+\int_{\OO} (r-\mathcal P r)_t \mathcal{P}e_r {\rm d}x\\
&\quad +\int_{\OO} (s-\mathcal P s)_t \mathcal{P}e_s {\rm d}x +\int_{\OO} (p-p_h)\mathcal{Q}e_p {\rm d}x+\int_{\OO} (q-q_h)\mathcal{Q}e_q {\rm d}x\\
&\quad +\mathcal B(r-\mathcal{P}r,p-\mathcal{Q}p,s-\mathcal{P}s,q-\mathcal{Q}q;~\mathcal{P}e_r,\mathcal{Q}e_p,\mathcal{P}e_s,\mathcal{Q}e_q)\\
&\quad +\mathcal B(\mathcal{P}e_r,\mathcal{Q}e_p,\mathcal{P}e_s,\mathcal{Q}e_q;~\mathcal{P}e_r,\mathcal{Q}e_p,\mathcal{P}e_s,\mathcal{Q}e_q)\\
&\quad =\mathcal H(r,s;~\mathcal{P}e_r,\mathcal{P}e_s)
-\mathcal H(r_h,s_h;~\mathcal{P}e_r,\mathcal{P}e_s).
\end{split}
\end{align*}

It follows from \eqref{pro} in Lemma \ref{lemma1} that
\begin{equation*}
  \mathcal B(r-\mathcal{P}r,p-\mathcal{Q}p,s-\mathcal{P}s,q-\mathcal{Q}q;~\mathcal{P}e_r,\mathcal{Q}e_p,\mathcal{P}e_s,\mathcal{Q}e_q)=0.
\end{equation*}
By the same argument as that used for the charge conservation law,
\begin{equation*}
\mathcal B(\mathcal{P}e_r,\mathcal{Q}e_p,\mathcal{P}e_s,\mathcal{Q}e_q;~\mathcal{P}e_r,\mathcal{Q}e_p,\mathcal{P}e_s,\mathcal{Q}e_q)=0.
\end{equation*}
The above two equations imply
\begin{align*}
\begin{split}
&\frac12\frac{d}{dt} \left[\int_{\OO}(\mathcal P e_r)^2{\rm d}x\right]
+\frac12\frac{d}{dt} \left[\int_{\OO}(\mathcal P e_s)^2{\rm d}x\right]
+\int_{\OO} (r-\mathcal P r)_t \mathcal{P}e_r {\rm d}x\\
&\quad +\int_{\OO} (s-\mathcal P s)_t \mathcal{P}e_s {\rm d}x +\int_{\OO} (p-p_h)\mathcal{Q}e_p {\rm d}x+\int_{\OO} (q-q_h)\mathcal{Q}e_q {\rm d}x\\[2mm]
&\quad =\mathcal H(r,s;~\mathcal{P}e_r,\mathcal{P}e_s)
-\mathcal H(r_h,s_h;~\mathcal{P}e_r,\mathcal{P}e_s).
\end{split}
\end{align*}

An application of the Cauchy-Schwarz inequality and Young inequality gives
\begin{align*}
\begin{split}
&\frac d{dt} \|\mathcal P e_r \|_{\mathbb L^2(\OO)}^2
+\frac d{dt} \|\mathcal P e_s \|_{\mathbb L^2(\OO)}^2+\|\mathcal Q e_p \|_{\mathbb L^2(\OO)}^2+\|\mathcal Q e_q \|_{\mathbb L^2(\OO)}^2\\[2mm]
&\le \|\mathcal P e_r\|_{\mathbb L^2(\OO)}^2+\|\mathcal P e_s\|_{\mathbb L^2(\OO)}^2+
\|(r-\mathcal P r)_t\|_{\mathbb L^2(\OO)}^2+
\|(s-\mathcal P s)_t\|_{\mathbb L^2(\OO)}^2\\[2mm]
&\quad +\|p-\mathcal Q p\|_{\mathbb L^2(\OO)}^2+
\|q-\mathcal Q q\|_{\mathbb L^2(\OO)}^2+2I,
\end{split}
\end{align*}
with
\begin{equation*}
\begin{split}
I&=|\mathcal H(r,s;~\mathcal{P}e_r,\mathcal{P}e_s)
-\mathcal H(r_h,s_h;~\mathcal{P}e_r,\mathcal{P}e_s)|\\
&\le\int_{\OO} |(r^2+s^2)r-(r_h^2+s_h^2)r_h|\cdot
|\mathcal{P}e_s|{\rm d}x\\
&\quad +\int_{\OO} |(r^2+s^2)s-(r_h^2+s_h^2)s_h|\cdot
|\mathcal{P}e_r|{\rm d}x
=:I_1+I_2.
\end{split}
\end{equation*}
Now, we give the error estimations of $I_1$ and $I_2$, respectively.

Since $u$ is bounded, it holds
\begin{align*}
&|(r^2+s^2)r-(r_h^2+s_h^2)r_h| \\
&=|(r^2+s^2+rr_h+r_h^2) e_r+(r_h s+r_h s_h) e_s| \\
&\le \|r^2+s^2+rr_h+r_h^2\|_{\mathbb L^\infty(\OO)}\cdot |e_r|
+\|r_h s+r_h s_h\|_{\mathbb L^\infty(\OO)}\cdot |e_s| \\
&\le C(1+\|u_h\|_{\mathbb L^\infty(\OO)})(|e_r|+|e_s|).
\end{align*}
Then H\"older and Cauchy-Schwarz inequalities yield
\begin{align*}
\begin{split}
I_1&\le C\Big(1+\|u_h\|_{\mathbb L^\infty(\OO)}\Big)\left(\int_{\OO} |e_r\mathcal P e_s| {\rm d}x
+\int_{\OO} |e_s\mathcal P e_s| {\rm d}x \right) \\
&\le C\Big(1+\|u_h\|_{\mathbb L^\infty(\OO)}\Big) 
\Big(\|s-\mathcal P s\|_{\mathbb L^2(\OO)}^2+\|r-\mathcal P r\|_{\mathbb L^2(\OO)}^2 \\
&\qquad \qquad \qquad \qquad \qquad +\|\mathcal P e_r\|_{\mathbb L^2(\OO)}^2+\|\mathcal P e_s\|_{\mathbb L^2(\OO)}^2\Big).
\end{split}
\end{align*}
In a similar manner, $I_2$ has analogous estimate. Thus
\begin{align*}
\begin{split}
I&\le C\Big(1+\|u_h\|_{\mathbb L^\infty(\OO)}\Big) 
\Big(\|s-\mathcal P s\|_{\mathbb L^2(\OO)}^2+\|r-\mathcal P r\|_{\mathbb L^2(\OO)}^2 \\
&\qquad \qquad \qquad \qquad \qquad +\|\mathcal P e_r\|_{\mathbb L^2(\OO)}^2+\|\mathcal P e_s\|_{\mathbb L^2(\OO)}^2\Big).
\end{split}
\end{align*}
Applying the interpolation properties \eqref{pro-est0} in Lemma \ref{pro-est}, we obtain\begin{align*}
& \frac d{dt} \|\mathcal P e_r \|_{\mathbb L^2(\OO)}^2
+\frac d{dt} \|\mathcal P e_s \|_{\mathbb L^2(\OO)}^2 \\
&\le C(1+\|u_h\|_{\mathbb L^\infty(\OO)})(\|\mathcal P e_r\|_{\mathbb L^2(\OO)}^2+\|\mathcal P e_s\|_{\mathbb L^2(\OO)}^2)\\
&\quad
+C(1+\|u_h\|_{\mathbb L^\infty(\OO)}) h^{2k+2} \Big( \|r\|_{\HH^{\bk+1} (\mathcal{I}_h)}^2
+\|s\|_{\HH^{\bk+1} (\mathcal{I}_h)}^2 \\
&\quad
+\|r_t\|_{\HH^{\bk+1} (\mathcal{I}_h)}^2
+\|s_t\|_{\HH^{\bk+1} (\mathcal{I}_h)}^2
+\|p\|_{\HH^{\bk+1} (\mathcal{I}_h)}^2
+\|q\|_{\HH^{\bk+1} (\mathcal{I}_h)}^2\Big).
\end{align*}
Finally, we conclude by Gronwall inequality and Lemma \ref{well} that
\begin{align}\label{87}
\begin{split}
&\sup_{t\in [0,T]}\|\mathcal P e_r \|_{\mathbb L^2(\OO)}
+\sup_{t\in [0,T]} \|\mathcal P e_s \|_{\mathbb L^2(\OO)} \\
&\le C \Big(1+\|u_h\|_{\mathbb L^\infty(\OO)}\Big) h^{\bk+1}
\Big( \sup_{t\in [0,T]}\|u\|_{\HH^{\bk+2} (\mathcal{I}_h)}
+\sup_{t\in [0,T]}\|u_t\|_{\HH^{\bk+1} (\mathcal{I}_h)} \Big)\\
&\le C\Big(1+\|u_h\|_{\mathbb L^\infty(\OO)}\Big) h^{\bk+1}.
\end{split}
\end{align}
Combining \eqref{pro-est0} and \eqref{87} and applying triangle inequality, we obtain the desired optimal error estimates \eqref{u-uh0}.
\quad \\

\begin{rk}
For linear Schr\"{o}dinger equation, the estimate \eqref{u-uh0} becomes
\begin{equation}
  \|u-u_h\|_{\mathbb L^2(\OO)}\le Ch^{\bk+1}.
\end{equation}
For nonlinear case, the above result is valid only when $u_h$ can be proved to be bounded.
It is not easy to be verified due to the discontinuity of the numerical solution.
It seems that some technical strategies are needed to derive the uniform boundedness of the numerical solution, and we will prove this claim in future work.
\end{rk}

\begin{rk}
We also note that the arguments in the present paper can be adapted to general NLS equation
\begin{align}\label{nls1}
\bi u_t+u_{xx}+f(|u|^2) u=0 \quad\text{in}\quad (0,T]\times \OO
\end{align}
with sufficiently smooth, bounded together with its first derivative, real-valued function $f$.
The main step in the proof of the optimal convergence rate (Theorem \ref{u-uh}) is the estimation of $I_1$.
In this case, set $g(r,s)=f(r^2+s^2)r$, then $g_r'=2f'(r^2+s^2)r^2+f(r^2+s^2)$ and
$g_s'=2f'(r^2+s^2)s^2+f(r^2+s^2)$.
Since $f$, $f'$ are assumed to be bounded, under appropriate assumptions on $u_0$, the solutions $u$ and $u_h$ of Eq. \eqref{nls} and Eq. \eqref{nls1}, respectively, can be shown to be bounded similarly to Lemma \ref{well}.
Applying Taylor expansion and H\"older and Cauchy-Schwarz inequalities yield
\begin{align*}
\begin{split}
I_1&\leq \int_{\OO} |(g_1)_{r}^{'}(\xi,\eta)||e_r\mathcal P e_s| {\rm d}x
+\int_{\OO} |(g_1)_{s}^{'}(\xi,\eta)||e_s\mathcal P e_s| {\rm d}x \\
&\leq C_1 \Big(1+\|u_h\|_{\mathbb L^\infty(\OO)}\Big)
\Big(\|s-\mathcal P s\|_{\mathbb L^{2}(\OO)}^2+\|r-\mathcal P r\|_{\mathbb L^{2}(\OO)}^2 \\
&\qquad \qquad \qquad \qquad \qquad +\|\mathcal P e_r\|_{\mathbb L^{2}(\OO)}^2+\|\mathcal P e_s\|_{\mathbb L^{2}(\OO)}^2\Big),
\end{split}
\end{align*}
where $\xi$ and $\eta$ lie between $r$ and $r_h$ and between $s$ and $s_h$, respectively.
\end{rk}

\section{Numerical Experiments}
\label{sec-4}

In this section, we will present some detailed numerical investigations of the CLDG method \eqref{ldg1}--\eqref{flux1} to the following NLS equation
\begin{align}\label{num-exp}
\bi u_t+u_{xx}+2|u|^2u=0.
\end{align}
Time discretization is by the implicit midpoint scheme. In particular, we will focus on the charge conservation law and accuracy of the method.

\subsection{The evolution of single soliton}

Consider Eq. \eqref{num-exp} with a single soliton solution
\begin{align*}
u(t,x)={\rm sech}(x+x_0-4t) e^{2\bi (x+x_0-3t/2)}.
\end{align*}

In the following experiments, we take the temporal step-size $\tau=0.001$, the spatial
meshgrid-size $h=0.5$, and the time interval $[0,~5]$, the numerical spatial domain $\OO=[-25,~25]$ with the periodic boundary condition.

\begin{figure}[th!]
\begin{minipage}[t]{0.5\linewidth}
\centering
\includegraphics[height=6cm,width=6cm]{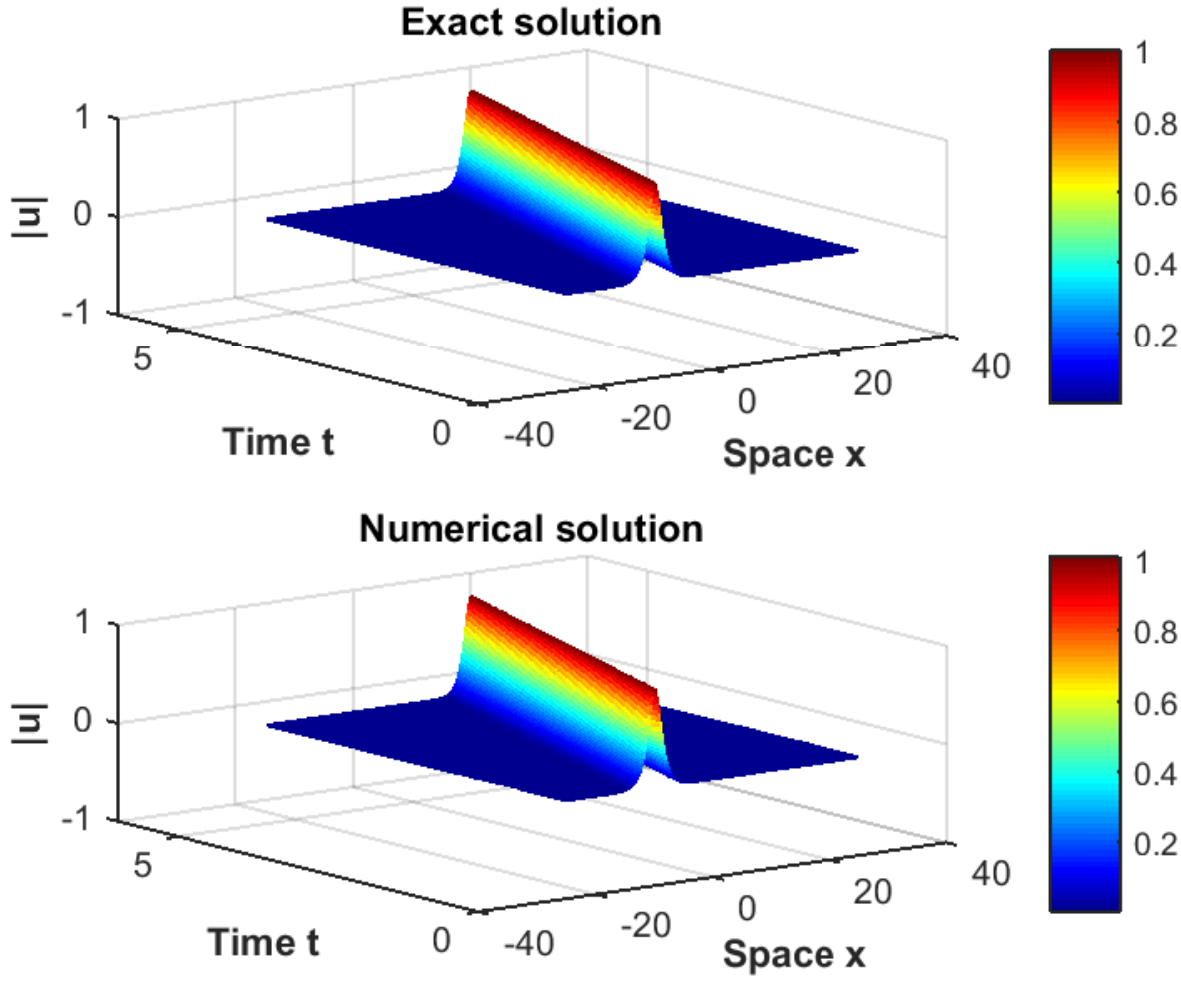}\\
\end{minipage}%
\begin{minipage}[t]{0.5\linewidth}
\centering
\includegraphics[height=6cm,width=6cm]{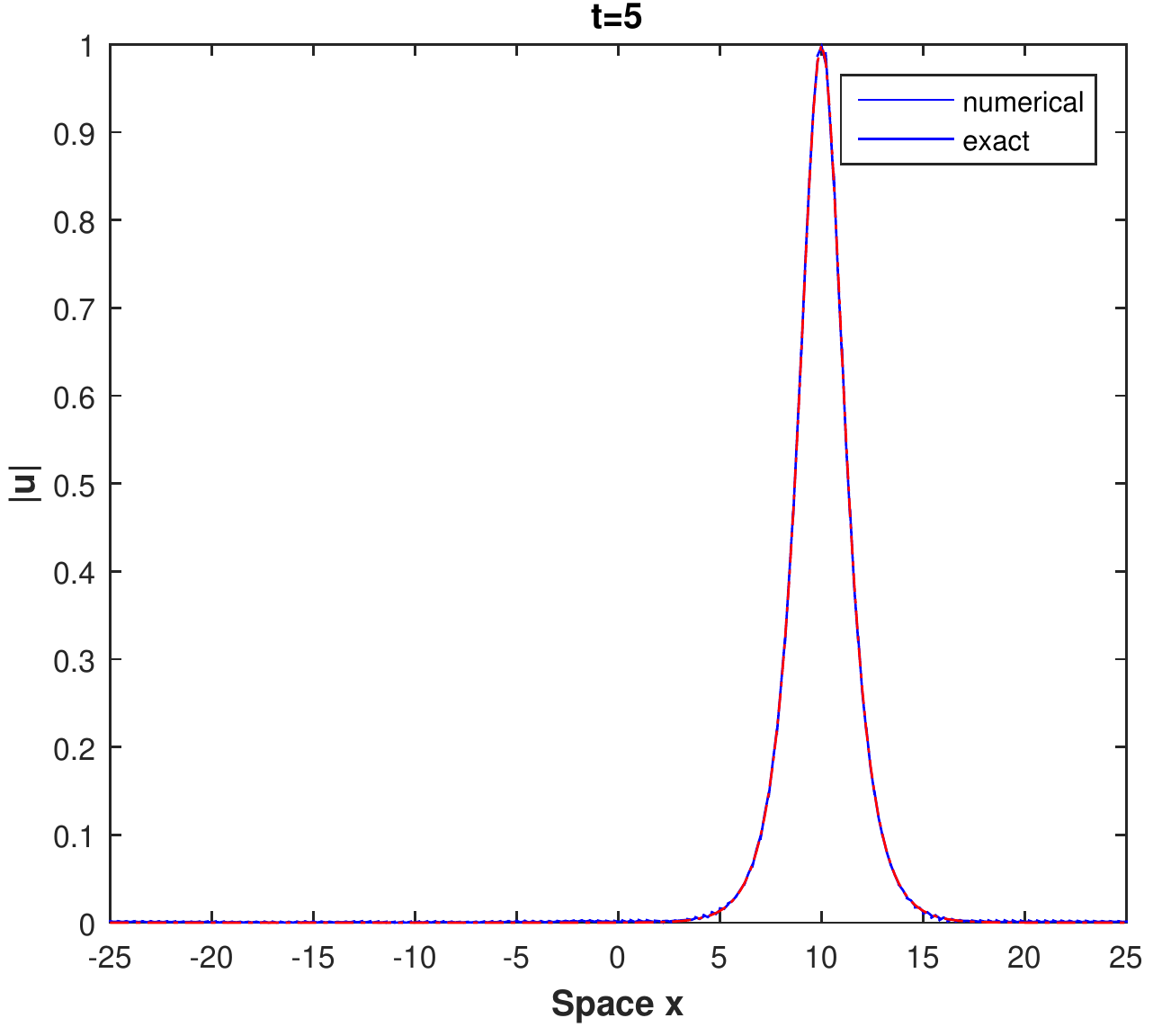}\\
\end{minipage}
\caption{Comparison of the analytic solution and the numerical solution with $x_0=10,\theta=1,h=0.5,\tau=0.001$ and $\bk=2$.}
\label{curl_single}
\end{figure}

The intensity profiles of the exact solution and numerical solution are shown in Fig. \ref{curl_single}. We observe a very good behavior of our method and a good agreement with the theoretical solution.

\begin{figure}[th!]
\centering
\includegraphics[height=7cm,width=8cm]{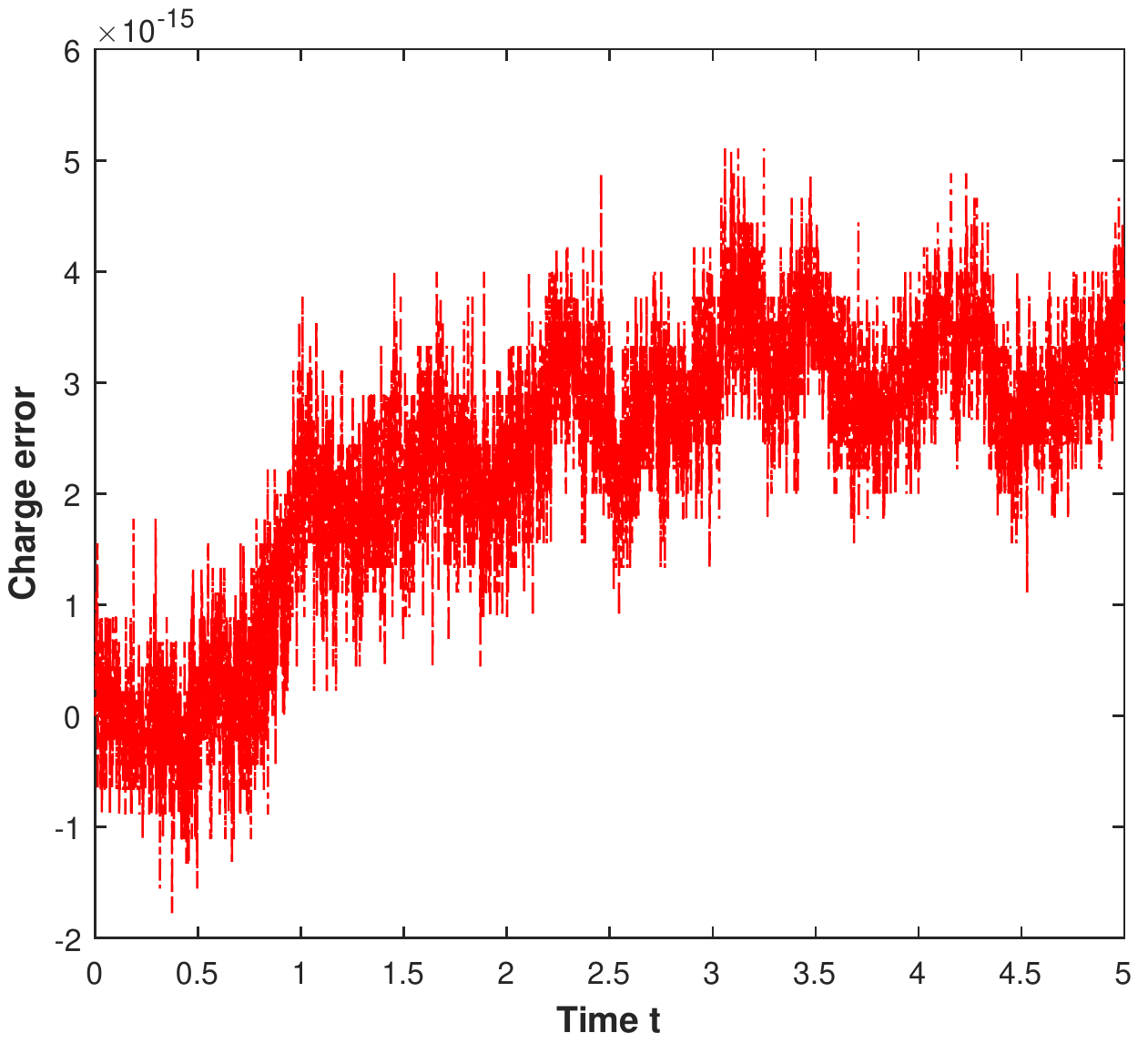}
\caption{The global error of the discrete charge conservation law of single soliton with $x_0=10,\theta=1,h=0.5,\tau=0.001$ and $\bk=2$.}
\label{charge_single}
\end{figure}

As is stated in Theorem 3.1, the LDG method \eqref{ldg0}-\eqref{flux} could preserve the discrete
charge conservation law exactly. We consider this phenomenon numerically in
Fig. \ref{charge_single}, where the figure shows the global error. We can see that the global residual of the
discrete charge conservation law reaches the magnitude of $10^{-15}$. Thus, we observe a good agreement with the theoretical result.

\subsection{The Interaction of Double Soliton}
In this experiment, we show the double soliton collision of Eq. \eqref{num-exp} with the initial condition
\begin{align*}
  u_0(x)={\rm sech}(x-x_1)e^{2\bi c_1(x-x_1)}
  +{\rm sech}(x-x_2)e^{2\bi c_2(x-x_2)}.
\end{align*}

\begin{figure}[th!]
\centering
\includegraphics[height=7cm,width=8cm]{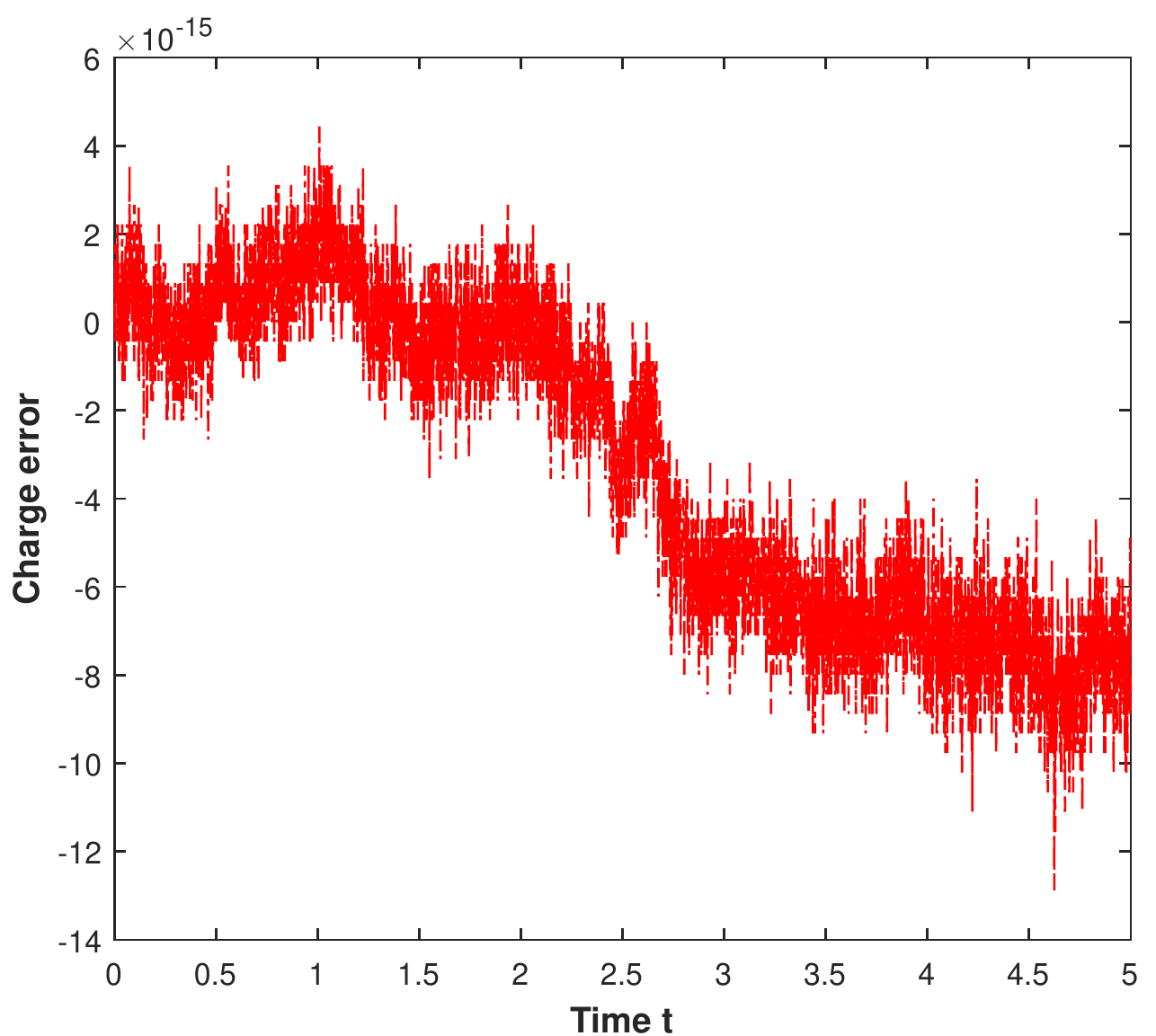}
\caption{The global error of the discrete charge conservation law of double soliton with $c_1=1$, $c_2=-1$, $x_1=-10$, $x_2=10,\theta=1,h=0.5,\tau=0.001$ and $\bk=2$.}\label{charge_double}
\end{figure}

The global error of the discrete charge conservation law is shown in Fig. \ref{charge_double}. Again, we observe the phenomena which agrees with the theoretical result.

\begin{figure}[th!]
\begin{center}
  \includegraphics[height=8cm,width=10cm]{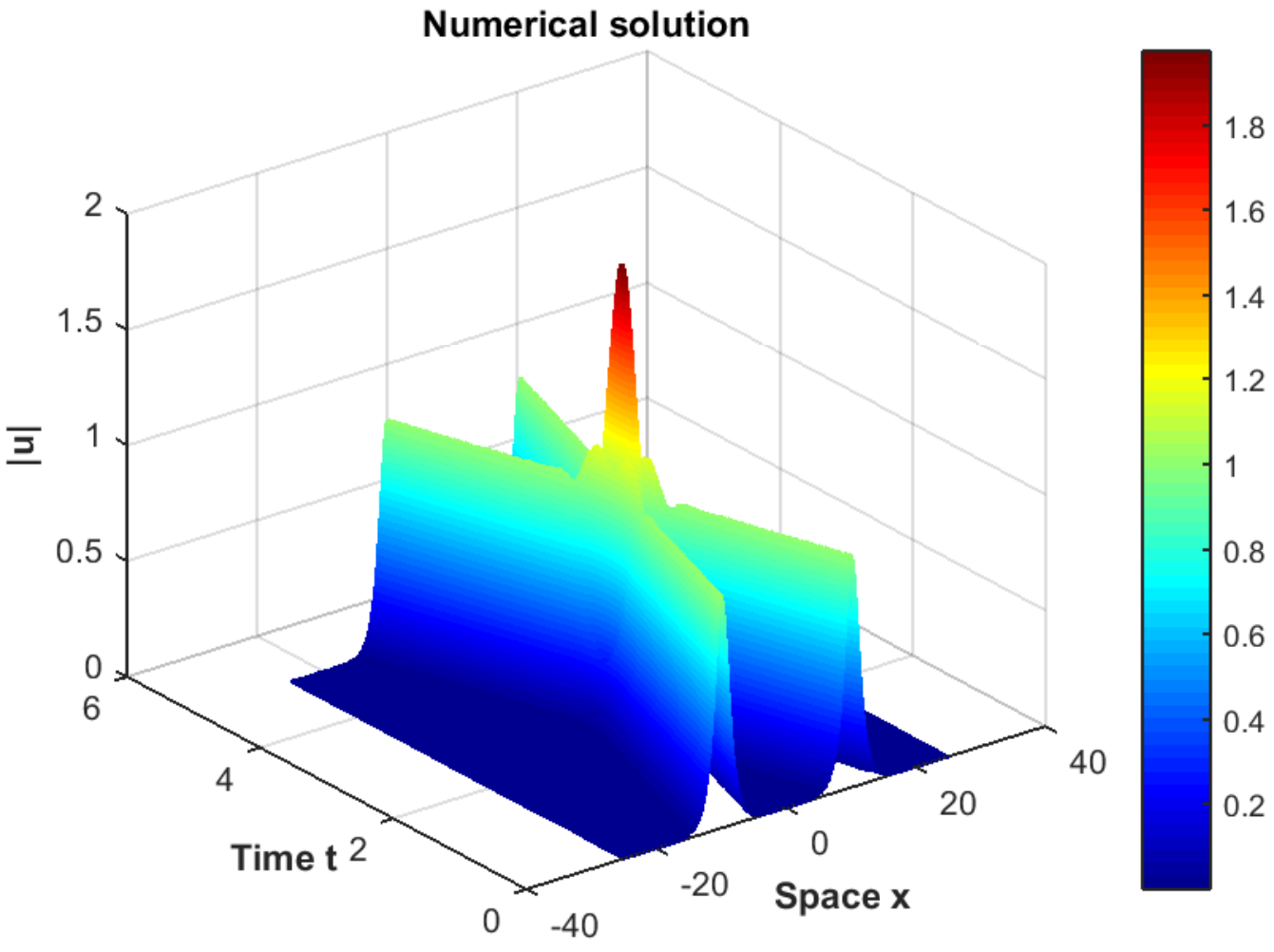}
  \caption{Collision of two solitons whose parameters are $c_1=1,c_2=-1,x_1=-10,x_2=10,\theta=1,h=0.2,\tau=0.001$ and $\bk=2$.}\label{curl_solution_double}
  \end{center}
\end{figure}

\begin{figure}[th!]
\begin{center}
  \includegraphics[height=12cm,width=12cm]{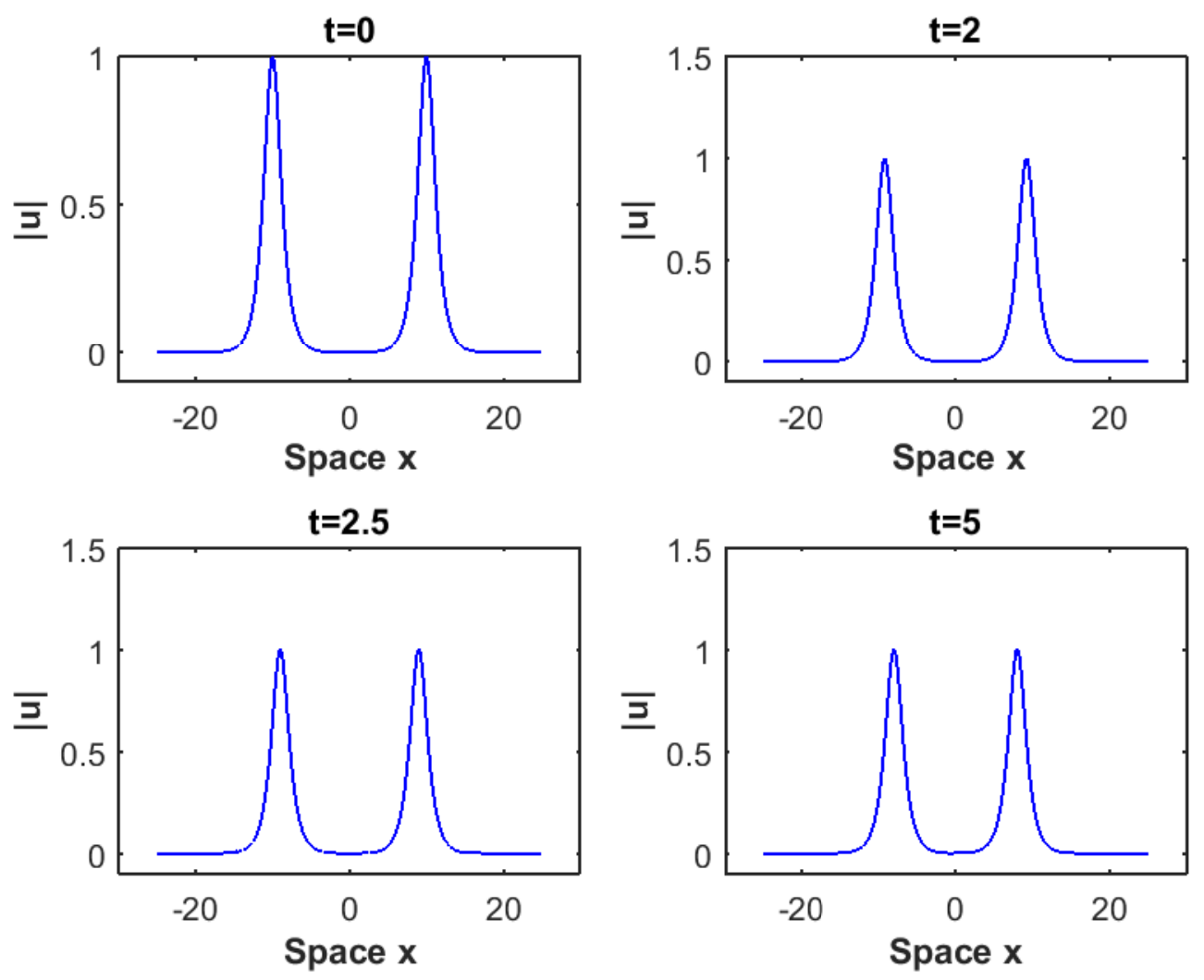}
  \caption{Profiles at time $t=0,~2,~2.5$ and $5$ of Fig. \ref{curl_solution_double}.} \label{interaction_double}
  \end{center}
\end{figure}

The evolution on the interval $[0; 5]$ is shown in Fig. \ref{curl_solution_double} and the profiles at different instants in Fig. \ref{interaction_double}. We observe that the interaction is elastic and the two waves emerges without any changes in their shapes and they conserve the energy almost exactly.

\subsection{The Birth of Mobile Soliton}
In this experiment, we show the birth of soliton using a square well initial condition
\begin{align*}
u_0(x)=Ae^{-x^2+2\bi x}.
\end{align*}

\begin{figure}[th!]
\begin{center}
  \includegraphics[height=8cm,width=10cm]{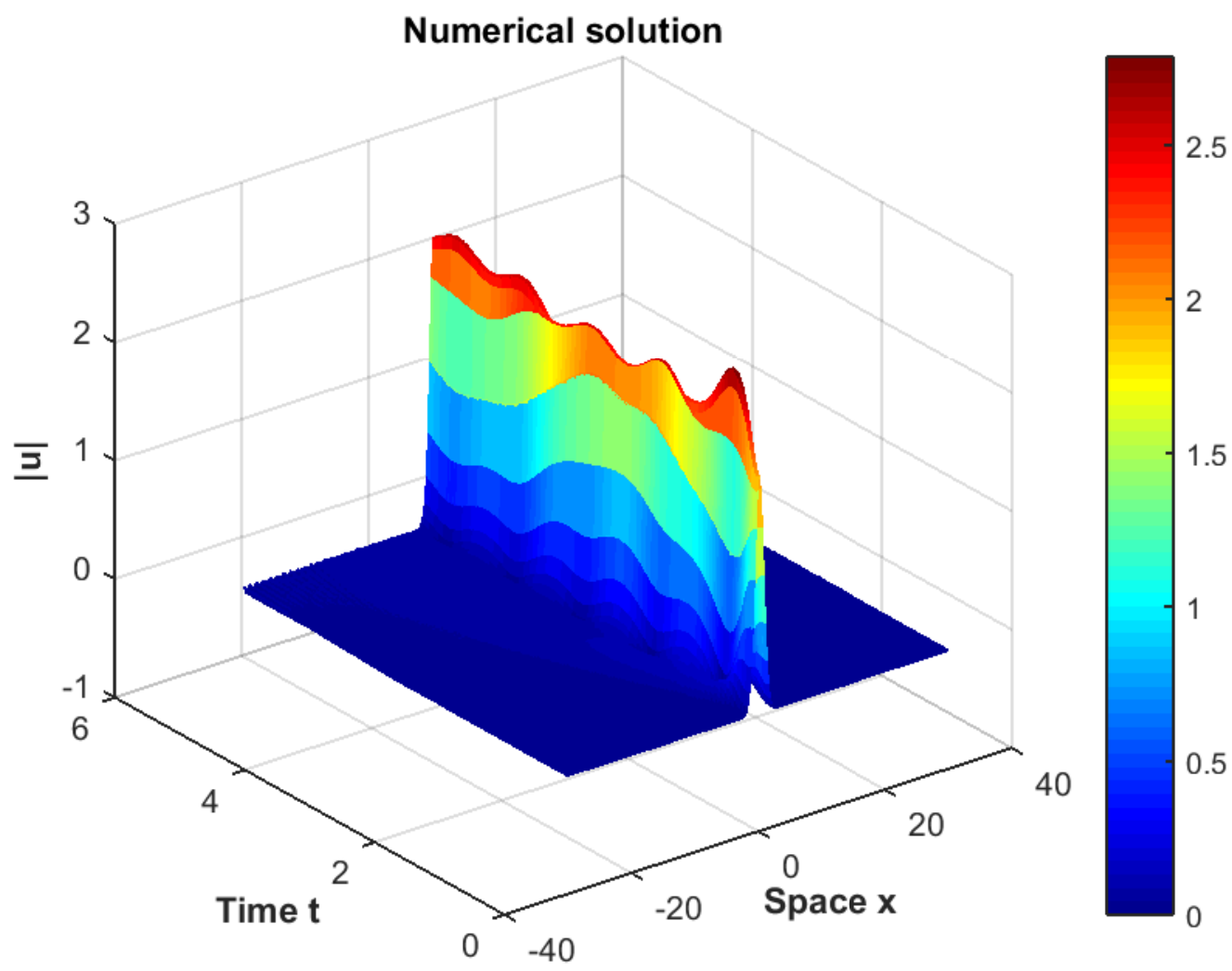}
  \caption{The birth of mobile soliton. $A=2,\theta=1,h=0.2,\tau=0.001$ and $\bk=2$. Periodic boundary condition in $[-30,30]$.}\label{curl_solution_square}
  \end{center}
\end{figure}

\begin{figure}[th!]
\begin{center}
  \includegraphics[height=12cm,width=12cm]{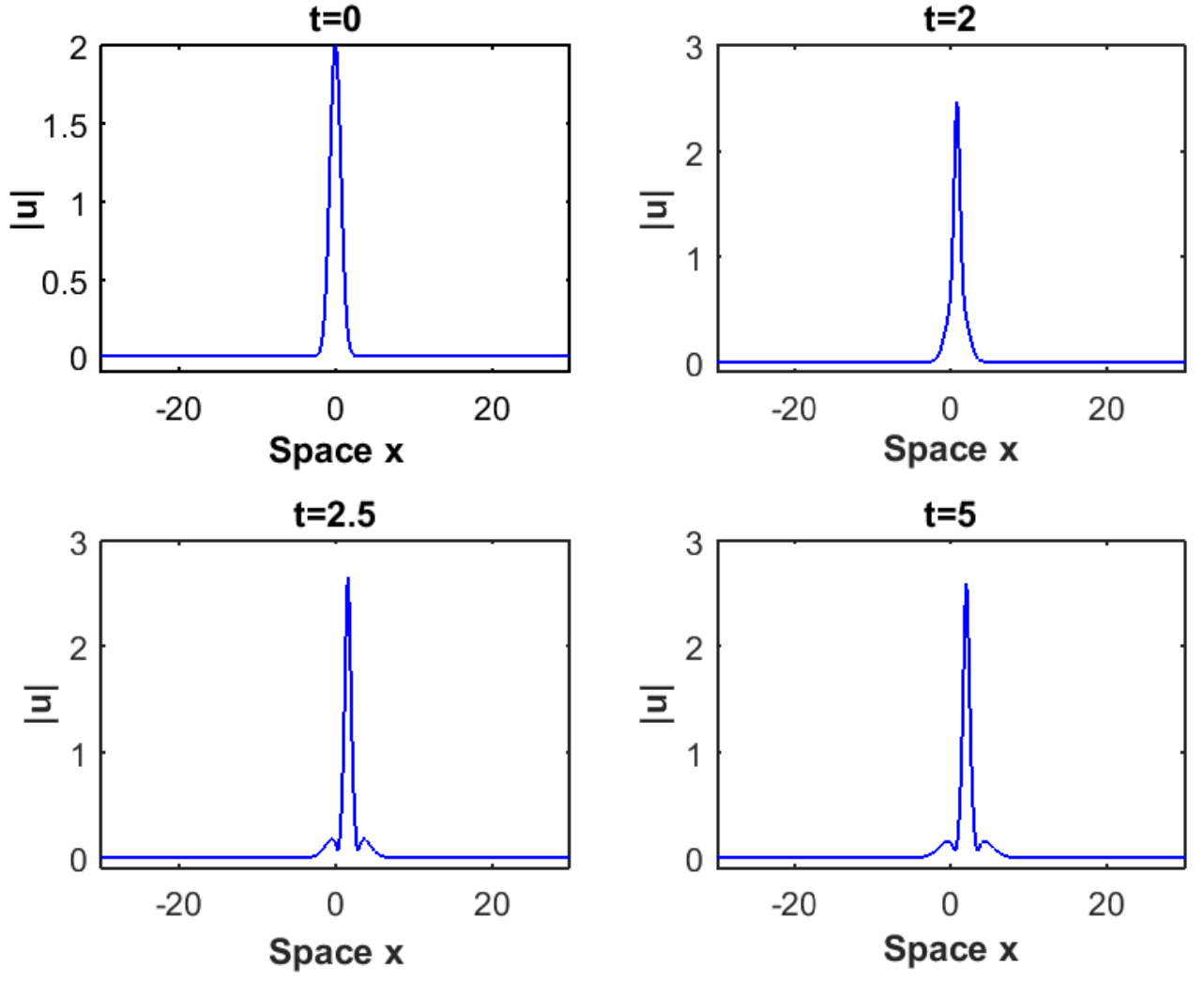}
  \caption{Profiles at time $t=0,~2,~2.5$ and $5$ of Fig. \ref{curl_solution_square}.}\label{birth_soliton}
  \end{center}
\end{figure}

In Figs. \ref{curl_solution_square}-\ref{birth_soliton}, the mobile soliton is observed.

\subsection{Optimal Convergence Estimates}
We show an accuracy test for Eq. \eqref{num-exp} with the soliton solution
\begin{align*}
u(t,x)={\rm sech}(x+x_0-4t)e^{2\bi (x+x_0-3t/2)}.
\end{align*}
We take the temporal step-size $\tau=0.00001$, and the time interval $[0,1]$, the numerical spatial domain $\OO=[-30,30]$ with the periodic boundary condition.

Table \ref{tab} lists the $\mathbb L^2$-errors and their numerical orders with different values of $\theta$ at $T=1$. From the table we conclude that, for all values of $\theta\in [0,1]$, one can always observe $(\bk+1)$-th order of accuracy in  $\mathbb L^2$-norm.

\begin{table}[th!]
\begin{minipage}[t]{0.5\linewidth}
\centering
\begin{tabular}{ccccc}
\hline   & $N$&$\mathbb L^2$-error ~~& Order~~\\[2mm]
\hline
                &60 &3.54E-1& -  \\[2mm]
                &120&7.20E-2&2.70\\[2mm]
    $\theta=0.4$&240&1.36E-2&3.02\\[2mm]
                &480&2.43E-3&3.04\\[2mm]

\hline
                &60 &2.31E-1& -  \\[2mm]
                &120&2.33E-2&3.73\\[2mm]
    $\theta=0.5$&240&3.89E-3&3.45\\[2mm]
                &480&6.91E-4&2.93\\[2mm]

\hline
                &60 &2.30E-1& -  \\[2mm]
                &120&2.89E-2&3.24\\[2mm]
    $\theta=1$&240&4.54E-3&3.68\\[2mm]
                &480&7.51E-4&2.94\\[2mm]
\hline
\end{tabular}
\end{minipage}%
\begin{minipage}[t]{0.5\linewidth}
\centering
\begin{tabular}{cccc}
\hline   & $N$&$\mathbb L^2$-error ~~& Order~~\\[2mm]
\hline
                &60 &8.72E-2& -  \\[2mm]
                &120&6.58E-3&3.94\\[2mm]
    $\theta=0.4$&240&5.64E-4&3.90\\[2mm]
                &480&4.82E-5&3.91\\[2mm]
\hline
                &60 &1.94E-2& -  \\[2mm]
                &120&2.18E-3&3.50\\[2mm]
    $\theta=0.5$&240&2.02E-4&3.87\\[2mm]
                &480&1.84E-5&3.80\\[2mm]

\hline
                &60 &2.91E-2& -  \\[2mm]
                &120&2.99E-3&3.84\\[2mm]
    $\theta=1$&240&2.64E-4&3.80\\[2mm]
                &480&2.67E-5&3.97\\[2mm]
\hline
\end{tabular}
\end{minipage}
\caption{Accuracy test for Eq. \eqref{nls} using
$\mathbb{P}^2$ polynomials (left) and $\mathbb{P}^{3}$ polynomials (right) with different $\theta$.}\label{tab}
\end{table}

\section{Conclusion}
\label{sec-5}

In this paper, we develop an LDG method to solve the one-dimensional nonlinear Schr\"{o}dinger equation.
The charge conservation law is shown to be preserved for LDG method proposed in this paper. The CLDG method, when applied to NLS equation, is shown to have the optimal $(\bk+1)$-th order of accuracy for polynomial elements of degree $\bk$.
The numerical tests demonstrate both accuracy and capacity of the method.

\section*{Acknowledgments}

The authors gratefully thank the anonymous referees for valuable comments and suggestions in improving this paper.
This work was supported by National Natural Science Foundation of China (No. 11601032, No. 91630312, No. 91530118, No. 11471310 and No. 11290142).

\bibliography{bib}

\end{document}